\documentstyle[11pt,leqno]{article}
\newtheorem{theorem}{{\sc Theorem}}
\newcommand{\bt}{\begin{theorem}}
\newcommand{\et}{\end{theorem}}
\newcommand{\newsection}[1]{\setcounter{equation}{0} \setcounter{theorem}{0}
\section{#1}}

\newcommand{\NI}{\noindent}
\newcommand{\bea}{\begin{eqnarray}}
\newcommand{\eea}{\end{eqnarray}}

\def \spec#1 {\mathop{#1}}

\def \b #1 {\bf #1}

\newcommand {\CC}{\centerline}

\newcommand{\clf}{{\cal F}}

\newcommand{\cly}{{\cal Y}}
\newcommand{\clz}{{\cal Z}}

\newcommand{\raro}{\rightarrow}

\newcommand{\vsp}{\vskip 1em}

\newcommand{\be}{\begin{equation}}
\newcommand{\ee}{\end{equation}}
\newcommand{\ben}{\begin{eqnarray*}}
\newcommand{\een}{\end{eqnarray*}}

\begin{document}
\CC{\bf{NONPARAMETRIC ESTIMATION OF LINEAR MULTIPLIER FOR} }
\CC{\bf{PROCESSES DRIVEN BY MIXED FRACTIONAL BROWNIAN MOTION}}
\vsp
\CC{\bf {B.L.S. Prakasa Rao}}
\CC{\bf {CR RAO Advanced  Institute of Mathematics, Statistics}}
\CC{\bf {and Computer Science, Hyderabad, India}}
\vsp
\NI{\bf Abstract:} We study the problem of nonparametric estimation of linear multiplier function $\theta(t)$ for processes satisfying stochastic differential equations of the type
$$dX_t= \theta(t)X_t dt+ \epsilon\;d\tilde W_t^H, X_0=x_0, 0 \leq t \leq T$$
where $\{\tilde W_t^H, t\geq 0\}$ is a mixed fractional Brownian motion with known Hurst index $H$ and study the asymptotic behaviour of the estimator as $\epsilon \raro 0.$
\vsp
\NI{\bf Keywords :} Nonparametric estimation, Linear multiplier, Mixed Fractional Brownian motion.
\vsp
\NI{\bf Mathematics Subject Classisfication :} Primary 62M09, Secondary 62G05.
\newsection {Introduction}

Statistical inference for fractional diffusion type processes satisfying stochastic differential equations driven by fractional Brownian motion  have been studied earlier and a comprehensive survey of various methods is given in Mishura (2008) and Prakasa Rao (2010). There has been a recent interest to study similar problems for stochastic processes driven by a mixed fractional Brownian motion (mfBm). Existence and uniqueness for solutions of stochastic differential equations driven by a mfBm are investigated in Mishura and Shevchhenko (2012) and Shevchenko (2014) among others. Maximum likelihood estimation for estimation of drift parameter in a linear stochastic differential equations driven by a mfBm is investigated in Prakasa Rao (2018a). The method of instrumental variable estimation for such parametric models is investigated in Prakasa Rao (2017). Some applications of such models in finance are presented in Prakasa Rao (2015 a,b). For related work on parametric inference for processes driven by mfBm, see Marushkevych (2016), Rudomino-Dusyatska (2003), Song and Lin (2014), Mishra and Prakasa Rao (2017), Prakasa Rao (2009) and Miao (2010) among others. Nonparametric estimation of the trend coefficient in models governed by stochastic differential equations driven by a mixed fractional Brownian motion is investigated in Prakasa Rao (2018b).

We now discuss the problem of estimating the function $\theta(t), 0 \leq t \leq T$ (linear multiplier) based on the observations of a  process $\{X_t, 0\leq t \leq T\}$ satisfying the stochastic differential equation
$$
dX_t=\theta(t)X_t dt + \epsilon \; d\tilde W^H_t, X_0=x_0, 0 \leq t \leq T
$$
where $\{\tilde W_t^H, t \geq 0\}$ is mfBM and study the properties of the estimator as $\epsilon \rightarrow 0.$

\vsp
\newsection {Mixed fractional Brownian motion}
We will now summarize some properties of stochastic processes which are solutions of stochastic differential equations driven by a mixed fractional brownian motion.
\vsp
Let $(\Omega, \clf, (\clf_t), P) $ be a stochastic basis
satisfying the usual conditions. The natural filtration of a
stochastic process is understood as the $P$-completion of the
filtration generated by this process. Let $\{W_t, t \geq 0\}$ be a standard Wiener process and $W^H= \{W_t^H, t \geq 0 \}$ be an independent  normalized  fractional Brownian motion with Hurst parameter $H \in (0,1)$, that is, a Gaussian process with continuous sample paths such that $W_0^H=0, E(W_t^H)=0$ and
\be
E(W_s^H W_t^H)= \frac{1}{2}[s^{2H}+t^{2H}-|s-t|^{2H}], t \geq 0, s \geq 0.
\ee
Let
$$\tilde W_t^H= W_t+ W_t^H, t \geq 0.$$
The process $\{\tilde W_t^H, t \geq 0\}$ is called the mixed fractional Brownian motion with Hurst index $H.$ We assume here after that Hurst index $H$ is known. Following the results in Cheridito (2001), it is known that the process $\tilde W^H$ is a semimartingale in its own filtration if and only if either $H=1/2$ or $H \in (\frac{3}{4},1].$
\vsp
Let us consider a stochastic process $Y=\{Y_t, t \geq 0\}$ defined
by the stochastic integral equation \be Y_t= \int_0^t C(s) ds  +\tilde W_t^H, t \geq 0 \ee where the process $C=\{C(t), t
\geq 0\}$ is an $(\clf_t)$-adapted process. For convenience, we write
the above integral equation in the form of a stochastic
differential equation
\be
dY_t= C(t) dt + d\tilde W_t^H, t \geq 0
\ee
driven by the mixed fractional Brownian motion $\tilde W^H.$ Following the recent works by Cai et al. (2016) and Chigansky and Kleptsyna (2015), one can construct an integral transformation that transforms the mixed fractional Brownian motion $\tilde W^H$ into a martingale $M^H.$ Let $g_H(s,t)$ be the solution of the integro-differential equation
\be
g_H(s,t)+H \frac{d}{ds}\int_0^t g_H(r,t)|s-r|^{2H-1} sign(s-r)dr=1, 0<s<t.
\ee
Cai et al. (2016) proved that the process
\be
M_t^H= \int_0^tg_H(s,t)d\tilde W_s^H, t \geq 0
\ee
is a Gaussian martingale with quadratic variation
\be
<M^H>_t= \int_0^tg_H(s,t)ds, t \geq 0
\ee
Furthermore the natural filtration of the martingale $M^H$ coincides with that of the mixed fractional Brownian motion $\tilde W^H.$ Suppose that, for the martingale $M^H$ defined by the equation (2.5), the sample paths of the process $\{C(t), t \geq 0\}$ are smooth enough in the sense that the process
\be
Q_H(t)= \frac{d}{d<M^H>_t}\int_0^tg_H(s,t)C(s)ds, t \geq 0
\ee
is well defined. Define the process
\be
Z_t= \int_0^tg_H(s,t) dY_s, t \geq 0.
\ee
As a consequence of the results in Cai et al. (2016), it follows that the process $Z$ is a fundamental semimartingale associated with the process $Y$ in the following sense.
\vsp
\NI{\bf Theorem 2.1:} {\it Let $g_H(s,t)$ be the solution of the equation (2.4). Define the process $Z$ as given in the equation (2.8). Then the following relations hold.

\NI(i) The process $Z$ is a semimartingale with the decomposition
\be
Z_t= \int_0^tQ_H(t)d<M^H>_s + M^H_t, t \geq 0
\ee
where $M^H$ is the martingale defined by the equation (2.5).

\NI(ii) The process $Y$ admits the representation
\be
Y_t=\int_0^t\hat g_H(s,t)dZ_s, t \geq 0
\ee
where
\be
\hat g_H(s,t)= 1-\frac{d}{d<M^H>_s}\int_0^tg_H(r,s)dr.
\ee
\NI (iii) The natural filtrations $(\cly_t)$ and $(\clz_t)$ of the processes $Y$ and $Z$ respectively coincide.}
\vsp
Applying Corollary 2.9 in Cai et al. (2016), it follows that the probability measures $\mu_Y$ and $\mu_{\tilde W^H}$ generated by the processes $Y$ and $\tilde W^H$ on an interval $[0,T]$ are absolutely continuous   with respect to each other and the Radon-Nikodym derivative is given by
\be
\frac{d\mu_Y}{d\mu_{\tilde W^H}}(Y)= \exp[\int_0^TQ_H(s) dZ_s-\frac{1}{2}\int_0^T[Q_H(s)]^2 d<M^H>_s]
\ee
which is also the likelihood function based on the observation $\{Y_s,0\leq s \leq T.\}$ Since the filtrations generated by the processes $Y$ and $Z$ are the same, the information contained in the families of $\sigma$-algebras  $(\cly_t)$ and $(\clz_t)$  is the same and hence the problem of the estimation of the parameters involved  based on the observation $\{Y_s, 0\leq s \leq T\}$ and $\{Z_s, 0\leq s \leq T\}$ are equivalent.
\vsp
\newsection{Preliminaries}
Let $\tilde W^H=\{W^H_t, t \geq o\}$ be a mixed fractional Brownian motion with known Hurst parameter $H.$ Consider the problem of estimating the function $\theta(t), 0 \leq t \leq T$ (linear multiplier) from the observations $\{X_t,0\leq t \leq T\}$ of process satisfying the stochastic differential equation
\be
dX_t=\theta(t) X_t dt + \epsilon \;d\tilde W^H_t, X_0=x_0, 0 \leq t \leq T
\ee
and study the properties of the estimator as $\epsilon \rightarrow 0 .$

Consider the differential equation in the limiting system of (3.1), that is , for $\epsilon=0,$ given by
\be
dx_t=\theta(t) x_t dt, x_0, 0 \leq t \leq T.
\ee
Observe that
$$x_t=x_0 \exp \{ \int^t_0 \theta(s) ds).$$

\vsp
We assume that the following condition holds:
\vsp
\NI{$(A_1)$:} The trend coefficient $\theta(t)$ over the interval $[0,T]$ is bounded by a constant $L$.
\vsp
The condition $(A_1)$  will  ensure the existence and uniqueness of the solution of the equation (3.1).
\vsp
\noindent{\bf Lemma 3.1:} Let the condition $(A_1)$ hold and  $\{X_t, 0\leq t \leq T\}$ and $\{x_t, 0\leq t \leq T\}$ be the solutions of the equations (3.1) and (3.2) respectively. Then, with probability one,
\be
|X_t-x_t| < e^{L t} \epsilon |\tilde W^H_t|
\ee
and
\be
\sup_{0 \leq t \leq T} E(X_t-x_t)^2  \leq e^{2L T} \epsilon^2 (T^{2H}+T).
\ee
\vsp
\noindent{\bf Proof of (a):} Let $u_t=|X_t-x_t|.$ Then by $(A_1)$; we have,
 \bea
u_t & \leq &\int^t_0 |\theta(v) (X_v-x_v)| dv + \epsilon |\tilde W^H_t|\\\nonumber
 & \leq & L \int^t_0 u_v dv + \epsilon |\tilde W^H_t|.\\\nonumber
\eea
Applying the Gronwall's lemma (cf. Lemma 1.12, Kutoyants (1994), p. 26), it follows that
\be
u_t \leq \epsilon |\tilde W^H_t| e^{L t}.
\ee
\vsp
\noindent{\bf Proof of (b):} From the equation (3.3), we have
\bea
E(X_t-x_t)^2 & \leq & e^{2 L t} \epsilon^2 E(|\tilde W^H_t|)^2\\\nonumber
&= & e^{2 L t} \epsilon^2 (t^{2H}+t). \\\nonumber
\eea
Hence
\be
\sup_{0 \leq t \leq T} E (X_t-x_t)^2 \leq e^{2 L T} \epsilon^2 (T^{2H}+T).
\ee
\vsp
Define
\bea
Q^*_{H, \theta} (t) & = &\frac{d}{d<M^H>_t}\int^t_0 g_H(t,s) \theta(s) x(s) ds \\\nonumber
&= &\frac{d}{d<M^H>_t}\int^t_0 g_H(t,s) \theta(s)[ x_0 \exp (\int^s_0 \theta(u) du)] ds\\\nonumber
\eea
by using the equation (3.2). We assume that
\vsp
\noindent{$(A_2)$:} the function  $Q^*_{H, \theta}(t)$ is Lipschitz of order $\gamma$ for any fixed $\theta (.).$
\vsp
Instead of estimation of the function $\theta (t),$ we consider the equivalent problem of estimating the function $Q^*_{H, \theta} (t)$ defined via the equation (3.9). This can be justified by the observation that the process $\{X_t, 0 \leq t \leq T\}$ governed by the stochastic differential equations (3.1)and the corresponding related process $\{Z_t, 0 \leq t \leq T\},$ as defined by  (2.8)  have the same filtrations by the results in Cai et al. (2016).
\vsp
We estimate the  function $Q^*_{H, \theta} (t)$ by a kernel type estimator defined by
\bea
\widehat{Q}_{H, \theta} (t) & = &\frac{1}{h_\epsilon} \int^T_0 G \left( \frac{s-t}{h_\epsilon}\right) dZ_s \\\nonumber
&= &\frac{1}{h_\epsilon}\int^T_0 G \left( \frac{s-t}{h_\epsilon}\right)\left( Q_{H, \theta} ds + \epsilon \;dM^H_s\right)\\\nonumber
\eea
by using the equation (2.9) where $G(u)$ is a bounded function with finite support $[A,B]$  satisfying the condition
\vsp
\noindent{ $(A_3)$: }$ G(u)=0$ for $ u < A, u > B$  and $\int^B_A G(u) du =1$.

It is obvious that the following conditions are satisfied by the function $G(.)$ \\
(i)$\displaystyle \int^\infty_{-\infty} G^2(u) du < \infty, $ and \\
(ii) $\displaystyle \int^\infty_{-\infty} G(u) |u|^\gamma du < \infty$  for $\gamma > 0,$ \\
\vsp

Consider a normalizing function $h_\epsilon \rightarrow 0$ with $\epsilon^2 h^{-3/2}_\epsilon \rightarrow 0$ as $\epsilon \rightarrow 0$.\\
\vsp
\newsection {Main Results}
\vsp
\noindent{\bf Theorem 4.1 :} Suppose the conditions $(A_1), (A_2)$ and $\ (A_3)$ are satisfied. Then for any $0 \leq t \leq T$, the estimator $\widehat{Q}_{H, \theta}(t)$ is uniformly consistent, that is,
\be
\lim_{\epsilon \rightarrow 0}\sup_{0 \leq t \leq T} E(|\widehat{Q}_{H, \theta}(t)-Q^*_{H, \theta} (t)|^2)=0.
\ee
\noindent{\bf Proof:}
\vsp
From (2.9), we have,
\bea
\;\;\;\\\nonumber
E |\widehat{Q}_{H, \theta} (t)-Q^*_{H, \theta} (t)|^2 & =& E \left|\frac{1}{h_\epsilon} \int^T_0 G \left(\frac{s-t}{h_\epsilon} \right)(Q_{H, \theta} (s) ds + \epsilon dM^H_s)-Q^*_{H, \theta} (t) \right|^2\\\nonumber
& = & E | \frac{1}{h_\epsilon} \int^T_0 G \left(\frac{s-t}{h_\epsilon} \right) (Q_{H, \theta}(s) -Q^*_{H, \theta}(s)) ds\\\nonumber
& & \;\;\;\;  + \frac{1}{h_\epsilon} \int^T_0 G \left(\frac{s-t}{h_\epsilon} \right) (Q^*_{H, \theta}(s) -Q^*_{H, \theta}(t)) ds\\\nonumber
& & \;\;\;\; + \frac{\epsilon}{h_\epsilon}\int^T_0 G \left(\frac{s-t}{h_\epsilon} \right)dM^H_s|^2 \\\nonumber
& = & E [I_1+I_2+I_3]^2 \;\;(\mbox{denoting the three integrals as $I_1, I_2$ and $I_3$ respectively})\\\nonumber
& \leq & 3\;E(I^2_1) + 3 \;E (I^2_2) + 3 \;E( I^2_3).\\\nonumber
\eea
Now
\bea
3\; EI^2_1 & = &3\;E \left|\frac{1}{h_\epsilon} \int^T_0 G \left( \frac{s-t}{h_\epsilon} \right) (Q_{H, \theta}(t) -Q_{H, \theta}(s)) ds \right|^2 \\\nonumber
& \leq &\frac{3}{h^2_\epsilon}[\int^T_0 G^2 \left( \frac{s-t}{h_\epsilon} \right) ds][ E \int^T_0 (Q_{H, \theta}(s)-Q^*_{H, \theta}(s))^2 ds].\\\nonumber
\eea
Now
\bea
\;\;\;\\\nonumber
\lefteqn{E \int^T_0 (Q_{H, \theta}(s) - Q^*_{H, \theta}(s))^2 ds} \\\nonumber
&=&\int^T_0 E \left[\frac{d}{d<M^H>_s} \int^s_0 g_H(s, v) \theta(v) (X(v) -x(v)) dv\right]^2 ds \\\nonumber
& = &\gamma_H \int^T_0 E \left[\int^s_0 \frac{\partial g_H(s,v)}{\partial s} \theta(v) (X(v)-x(v)) dv\right]^2 \beta(s,H) ds\\\nonumber
& \leq &  \gamma_H \int^T_0 \beta(s,H)  \{ \int^s_0 \left ( \frac{\partial g_H(s,v)}{\partial s} \right)^2 \theta^2(v)dv \int^s_0 E(X(v)-x(u))^2 dv\}ds\\\nonumber
\eea
where the constant $\gamma_H$ and the function $\beta(s,H)$ depend on the quadratic variation of the martingale $M^H.$
Note that
$$
E(X_v-x_v)^2 \leq e^{2 L v} \epsilon^2( v^{2H}+v) \;\;\mbox{(by Lemma 3.1)}.
$$
Hence, from the equation (4.4) and $(A_3)$ we get that,
\bea
3EI^2_1  &\leq & C \frac{(B-A)}{h^2_\epsilon} \left\{ \int^\infty_{-\infty} G^2(u) du\right\} \epsilon^2 h_\epsilon \\\nonumber
& & \;\;\;\;\times \int^T_0 \beta(s,H) \left\{\int^s_0 e^{2 L v} (v^{2H}+v)  dv \right\}\left\{ \int^s_0 \left ( \frac{\partial g_H(s,v)}{\partial s} \right)^2 dv \right\} ds \\\nonumber
& \leq & \epsilon^2 h^{-1}_\epsilon C(T, L, H) \\\nonumber
\eea
and the last term tends to zero as $\epsilon \rightarrow 0.$
\vsp
In addition,
\bea
I^2_2&= &3 \left\{\frac{1}{h_\epsilon} \int^T_0 G \left( \frac{s-t}{h_\epsilon}\right) (Q^*(s) ds-Q^*(t)) ds\right\}^2 \\\nonumber
&= &3 \left\{\int^\infty_{-\infty} G(u) (Q^*(t+h_\epsilon u)-Q^*(t)) du \right\}^2 \;\;(\mbox{by}  (A_2))\\\nonumber
& \leq &C \left\{\int^\infty_{-\infty}G(u) |h_\epsilon u|^\gamma  du \right\}^2 \;\;(\mbox{by}  (A_2))\\\nonumber
& \leq &C h^{2\gamma}_\epsilon \left( \int^\infty_{-\infty}G(u) |u|^\gamma du\right)^2\\\nonumber
& \leq &C h^{2\gamma}_\epsilon \;\;\mbox{by}  (A_3))\\\nonumber
\eea
and the last term tends to zero as  $\epsilon \rightarrow 0.$ Furthermore
\bea
I^2_3&= &\frac{3 \epsilon^2}{h_{\epsilon}^2}E \left( \int^T_0 G \left( \frac{s-t}{h_\epsilon} \right) dM^H_s \right)^2\\\nonumber
& = &\frac{3 \epsilon^2}{h_{\epsilon}^2}\int^T_0 G^2 \left( \frac{s-t}{h_\epsilon} \right) d<M^H>_s \\\nonumber
& = &\frac{3 \epsilon^2}{h_{\epsilon}^2}\int^T_0 G^2 \left( \frac{s-t}{h_\epsilon} \right) \beta(s,H) ds\\\nonumber
& \leq &\frac{3 \epsilon^2}{h_{\epsilon}^2} \left\{\int^T_0 G^2  \left( \frac{s-t}{h_\epsilon} \right)ds \int^T_0 \beta^2(s,H)ds\right\}^\frac{1}{2}\\\nonumber
& \leq &C \frac{3 \epsilon^2}{h_{\epsilon}^2} \left\{h_\epsilon (\int^\infty_{-\infty} G^2(u)  \ du) D(H,T) \right\}^\frac{1}{2}\\\nonumber
& \leq &C(T,H) \epsilon^2 h^{-3/2}_\epsilon .\\\nonumber
\eea
for some constants $C, D(T,H)$ and $C(T,H)$ depending on $T$ and $H.$ The result follows from the equations (4.5), (4.6) and (4.7).
\vsp
\noindent{\bf Corollary 4.2:} Under the conditions $(A_1), (A_2)$ and $\ (A_3),$

$$\displaystyle \lim_{\epsilon \rightarrow 0} \sup_{|\theta (.)| \leq L, 0 \leq t \leq T}E \left\{ \widehat{Q}_{H, \theta}(t)-Q^*_{H, \theta}(t)\right
\}^2 \overline{\epsilon}^{\frac{8\gamma}{4\gamma+3}}< \infty.$$\\
\vsp
\noindent{\bf Proof:} From the inequalities derived in (4.5), (4.6) and (4.7), we get that
there exist positive constants $C_1,C_2$ and $C_3$ depending on $T$ and $H$ such that
\be
 \sup_{|\theta(.)| \leq L, 0 \leq t \leq T}E \left\{\widehat{Q}_{H, \theta}(t)-Q^*_{H, \theta} (t)\right\}^2 \leq C_1 \epsilon^2 h^{-1}_\epsilon + C_2 h^{2\gamma}_\epsilon +C_3 \epsilon^2 h^{-\frac{3}{2}}_\epsilon.
\ee
Let $h_\epsilon=\epsilon^\beta, 0<  \beta  < \frac{4}{3}. $ Then the condition $h^{2\gamma}_\epsilon =\epsilon^{2}h_\epsilon^{-3/2}$ leads to the choice
$\beta=\frac{4}{4 \gamma+3}$ and  we get an optimum bound in (4.8) and hence
\be
\lim_{\epsilon \rightarrow 0} \sup_{|\theta (.)| \leq L, 0 \leq t \leq T} E \left[ \widehat{Q}_{H, \theta}(t)-Q^*_{H, \theta}(t) \right]^2\epsilon^{-\frac{8 \gamma}{4 \gamma +3}} \leq C
\ee
for some positive constant $C$ which implies the result.
\vsp
\NI{\bf Acknowledgement:} This work was supported by the Indian National Science Academy (INSA) under the scheme  ``INSA Senior Scientist" at the CR RAO Advanced Institute of Mathematics, statistics and Computer Science, Hyderabad 500046, India.
\vsp
\NI{\bf References :}
\vsp
\begin{description}

\item Cai, C., Chigansky, P. and Kleptsyna, M. (2016) Mixed Gaussian processes, {\it Ann. Probab.}, {\bf 44}, 3032-3075.

\item Cheridito, P. (2001) Mixed fractional Brownian motion, {\it Bernoulli}, {\bf 7}, 913-934.

\item Chigansky, P. and Kleptsyna, M. (2015) Statistical analysis of the mixed fractional Ornstein-Uhlenbeck process, arXiv:1507.04194.

\item Kutoyants, Y.A.(1994) {\it Identification of Dynamical Systems with Small Noise}, kluwer, Dordrecht.

\item Marushkevych, Dmytro. (2016) Large deviations for drift parameter estimator of mixed fractional Ornstein-Uhlenbeck process, {\it Mod.  Stoch. Theory Appl.}, {\bf 3}, 107-117.

\item  Miao, Y. (2010) Minimum $L_1$-norm estimation for mixed Ornstein-Uhlenbeck type process, {\it Acta. Vietnam}, {\bf 35}, 379-386.

\item Mishura, Y. (2008){\it Stochastic Calculus for Fractional Brownian Motion and Related Processes}, Berlin:Springer.

\item Mishura, Y. and Shevchenko, G. (2012) Existence and uniqueness of the solution of stochastic differential equation involving Wiener process and fractional Brownian motion with Hurst index $H> 1/2,$ {\it Comput. Math. Appl.}, {\bf 64}, 3217-3227.

\item Mishra, M.N. and Prakasa Rao, B.L.S. (2017) Large deviation probabilities for maximum likelihood estimator and Bayes estimator of a parameter for mixed fractional Ornstein-Uhlenbeck type process, {\it Bull. Inform. and Cyber.}, {\bf 49}, 67-80.

\item Prakasa Rao, B.L.S. (2009) Estimation for stochastic differential equations driven by mixed fractional Brownian motion, {\it Calcutta Stat. Assoc. Bull.}, {\bf 61}, 143-153.

\item Prakasa Rao, B.L.S. (2010) {\it Statistical Inference for Fractional Diffusion Processes}, Wiley, Chichester.

\item Prakasa Rao, B.L.S. (2015a) Option pricing for processes driven by mixed fractional Brownian motion with superimposed jumps, {\it Probability in the Engineering and Information Sciences}, {\bf 29},  589-596.

\item Prakasa Rao, B.L.S. (2015b) Pricing geometric Asian power options under mixed fractional Brownian motion environment, {\it Physica A}, {\bf 446}, 92-99.

\item Prakasa Rao, B.L.S. (2017) Instrumental variable estimation for a linear stochastic differential equation driven by a mixed fractional Brownian motion, {\it Stochastic Anal. Appl.}, {\bf 35}, 943-953.

\item Prakasa Rao, B.L.S. (2018a) Parameter estimation for linear stochastic differential equations driven by  mixed fractional Brownian motion, {\it Stochastic Anal. Appl.}(to appear).

\item Prakasa Rao, B.L.S. (2018b) Nonparametric estimation of trend for stochastic differential equations driven by mixed fractional Brownian motion, Preprint.

\item Rudomino-Dusyatska, N. (2003) Properties of maximum likelihood estimates in diffusion and fractional Brownian models, {\it Theor. Probab. Math. Statist.}, {\bf 68}, 139-146.

\item Song, N. and Liu, Z. (2014) Parameter estimation for stochastic differential equations driven by mixed fractional Brownian motion, {\it  Abst. Appl. Anal.} : 2014 Article ID 942307, 6 pp.

\item Shevchenko, G. (2014) Mixed stochastic delay differential equations, {\it Theory Probab. Math. Statist.}, {\bf 89}, 181-195.

\end{description}

\end{document}